\documentclass[12pt]{article}
\usepackage[cp1251]{inputenc}
\usepackage[russian]{babel}
\usepackage{amssymb,latexsym,amsmath}
\usepackage{oldlfont}
\usepackage{amsfonts}
\usepackage{enumerate}
\usepackage {graphicx,color}
\usepackage{wrapfig}
\usepackage{amsfonts,amstext,amsmath,float}

\oddsidemargin\evensidemargin
\usepackage{array,longtable}
\newcommand{\D}{\displaystyle}

\definecolor{MyRed}{rgb}{1,0,0}
\definecolor{MyGreen}{rgb}{0,1,0}
\definecolor{MyBlue}{rgb}{0,0,1}

\begin{document}

\begin{center}
{\Large\bf A class  of iterative methods for solving nonlinear operator equations}
\end{center}

\begin{center}
{\sc O.N. Evkhuta, P.P. Zabre\u{\i}ko}
\end{center}

\centerline{\bf Summary}

\begin{center}
\parbox{15.0cm}{{\small \hspace{0.5cm} The article deals with gradient-like iterative methods for solving nonlinear operator equations on Hilbert and Banach spaces. The authors formulate a general principle of studying such methods. This principle allows  to formulate simple conditions of convergence of the method under consideration, to estimate the rate of this convergence and to give effective {\small\sc apriori} and {\small\sc aposteriori} error estimates  in terms of a scalar function that is constructed on the base of estimates for properties of invertibility and smoothness of linearizations of the left-hand side of the equations under study. The principle is applicable for analysis of such classical methods as method of minimal residuals, method of steepest descent, method of minimal errors and others. The main results are obtained for operator equations on Hilbert spaces and Banach spaces with a special property, that is called Bynum property.}}
\end{center}

\vspace{0.3cm}

The method of steepest descent \cite{KA}, method of minimal residuals \cite{KK}, method of minimal errors \cite{FR1,FR2} and some other gradient-like methods are simple and effective tools of studying and  solving linear operator equations on Hilbert spaces. In the articles of M. Altman (\cite{9,10,11}), L.V. Kivistik (\cite{K2}), V.M. Fridman (\cite{FR3})  some modifications of these methods for nonlinear operator equations  were investigated. An account of these methods for linear and nonlinear operator equations was presented
in the book \cite{KVZRS}.
All these methods are based on a simple idea that the passage from the iteration $x_n$ to the iteration $x_{n+1}$, for each step, is realized in such a way that the value of residual at $x_{n+1}$ is  less than the value residual at $x_n$. It turns out  \cite{5,6,2,1,3} that, under the standard conditions of convergence for the classical methods of steepest descent, minimal residuals, minimal errors, and others, the analysis of these methods is reduced to the construction of a scalar function $d(\cdot)$ (sometimes it is called {\it the relaxation function}) such that the inequalities $\|f(x_{n+1})\| < d(\|f(x_n)\|)$, $n = 0,1,2,\ldots$, hold. Moreover, the graph of this function allows us to give a sufficiently full characterization of the method in use:  conditions of its convergence and rate of this convergence, {\small\sc apriori} and {\small\sc aposteriori} error estimates and so on. In addition, these conditions are usually more general than the classic ones and the corresponding error estimates are more proximate. At last, the approach, based on the analysis of function $d(\cdot)$ can be extended onto operator equations on Banach spaces.

In this paper we present the analysis of  the approach described above within the frames of  investigation of gradient-like iterative methods and obtain the corresponding convergence conditions of these methods, estimates for the rate of this convergence and {\small\sc apriori} and {\small\sc aposteriori} error estimates of the corresponding approximations.

\vspace{0.5cm}

{\bf 1.} Let $X$ be a Hilbert space and let $B(x_0,R) \ (\subset X)$ denote the ball with the center $x_0$ and of the radius $R$.  We are interested in the solvability of the equation
 \begin{equation}\label{equation}
 f(x) = 0
 \end{equation}
in the ball $B(x_0,R)$ and in the convergence of the approximations
 \begin{equation}\label{approximations}
 x_{n+1} = x_n - \Lambda(x_n,f(x_n))T(x_n)f(x_n) \quad (n =  1,2,\ldots)
 \end{equation}
where
 \begin{equation*}
 f: \ B(x_0,R) \ (\subset X) \to X,
 \end{equation*}
 \begin{equation*}
 T: \ B(x_0,R) \ (\subset X) \to L(X),
 \end{equation*}
 \begin{equation*}
 \Lambda: \ B(x_0,R) \times X \ (\subset X \times X) \to \mathbb{R}.
 \end{equation*}
We assume that $f$ is Fr\'echet-differentiable at every point of $B(x_0,R)$ and $T, \Lambda, f$ satisfy the  conditions
 \begin{equation}\label{CondLambda}
 \Lambda(x,f(x)) \le \lambda(r), \quad \|T(x)\| \le \theta(r) \qquad (\|x - x_0\| \le r, \  0 \le r \le R)
 \end{equation}
 \begin{equation}\label{CondLambdaT}
 {\|f(x)\|}^{-1} \|f(x) - \Lambda(x,f(x))f'(x)T(x)f(x)\| \le \mu(r) \qquad (\|x - x_0\| \le r, \  0 \le r \le R)
 \end{equation}
 \begin{equation}\label{Condomega}
 \begin{array}{c}\|f'(x_1) - f'(x_2)\| \le \omega(r,\|x_1 - x_2\|)  \\[4pt] (\|x_1 - x_0\|, \|x_2 - x_0\| \le r, \  0 \le r \le R).\end{array}
 \end{equation}

Let us introduce the following functions:
 \begin{equation}\label{Fun d}
 d(r,\phi) = \mu(r)\phi + \Omega(r,\lambda(r)\theta(r)\phi) \qquad \biggl(\Omega(r,t) = \int\limits_0^t \omega(r,\tau) \, d\tau\biggr),
 \end{equation}
 \begin{equation}\label{Fun d^n}
 d^{(0)}(r,\phi) = \phi, \qquad d^{(n+1)}(r,\phi) =  d^{(n)}(r,d(r,\phi))\quad (n =1,2,\ldots).
 \end{equation}
These function are defined on the interval $[0,R]$, the function $d(r,\phi)$ is convex for small positive $\phi$; its graph is presented at the following picture:

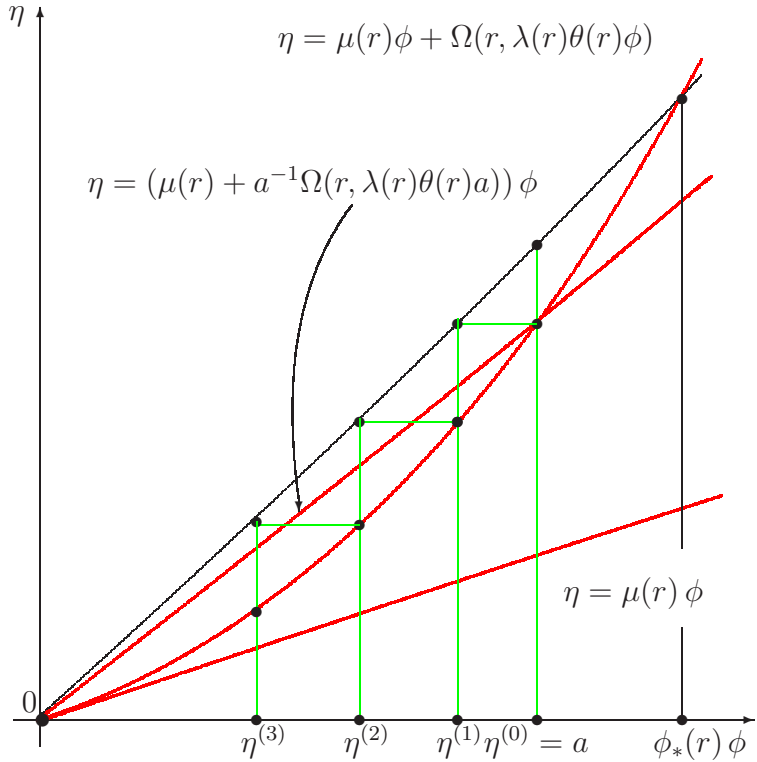
\begin{figure}[H]
\begin{center}
\begin{picture}(300,300)
\put(10,-5){\vector(0,1){280}} \put(0,5){\vector(1,0){280}}
\qbezier(8,5)(150,133)(260,249)
\put(270,-7){\text{$\phi$}}
\put(-2,270){\text{$\eta$}} \put(3,8){\text{$0$}}
{\color{MyRed}\thicklines\qbezier(9,5)(150,50)(260,255)}
{\color{MyRed}\thicklines\qbezier(5,5)(150,118)(260,211)}
{\color{MyRed}\thicklines\qbezier(3,5)(150,55)(260,90)}
\thinlines\put(245,5){\line(0,1){35}}
\thinlines\put(245,70){\line(0,1){170}}
\put(245,5){\circle*{4}} \put(245,240){\circle*{4}}
\put(234,-7){$\phi_*(r)$}
\put(92,260){$\eta =  \mu(r)\phi + \Omega(r,\lambda(r)\theta(r)\phi)$}
{\color{MyGreen}\put(190,5){\line(0,1){180}}}
\put(190,5){\circle*{4}} \put(190,155){\circle*{4}} \put(190,185){\circle*{4}}
\put(170,-7){$\eta^{(0)} = a$}
\put(200,50){$\eta = \mu(r)\,\phi$}
{\color{MyGreen}\put(190,155){\line(-1,0){30}}} \put(160,155){\circle*{4}}
{\color{MyGreen}\put(160,5){\line(0,1){152}}} \put(160,5){\circle*{4}}
\put(160,118){\circle*{4}} \put(152,-7){$\eta^{(1)}$}
{\color{MyGreen}\put(160,118){\line(-1,0){37}}} \put(160,118){\circle*{4}}
\put(123,118){\circle*{4}} \put(117,-7){$\eta^{(2)}$}
{\color{MyGreen}\put(123,5){\line(0,1){114}}} \put(123,79){\circle*{4}}
\put(123,5){\circle*{4}} \put(3,5){\circle*{5}}
{\color{MyGreen}\put(123,79){\line(-1,0){39}}} \put(84,80){\circle*{4}}
{\color{MyGreen}\put(84,5){\line(0,1){75}}} \put(84,5){\circle*{4}}
\put(78,-7){$\eta^{(3)}$} \put(84,46){\circle*{4}}
\put(20,205){$\eta = (\mu(r) + a^{-1}\Omega(r,\lambda(r)\theta(r)a))\,\phi$}
\qbezier(120,200)(90,160)(100,88)
\put(100,88){\vector(0,-2){4}}
\end{picture}
\end{center}
\end{figure}

Let $\phi_*(r)$ be the smallest root at $[0,R]$ of the equation
 \begin{equation}\label{Equation for psi}
 \phi = d(r,\phi)
 \end{equation}
(we assume, that this root does exist). Put
 \begin{equation}\label{Fun w}
 w(r,\phi) = \sum_{n=0}^\infty d^{(n)}(r,\phi) \qquad (0 \le \phi <  \phi_*(r)),
 \end{equation}
This function is defined on the interval $(0,\phi_*(r))$ and possesses the following important property:
 \begin{equation}\label{Equation for w}
 w(r,\phi) = \phi + w(r,d(r,\phi)) \qquad (0 \le \phi < \phi_*(r))
 \end{equation}
Moreover, the following inequality
 \begin{equation*}
 w(r,\phi) \le \frac{\phi^2}{\phi - d(r,\phi)} \qquad (0 \le \phi <  \phi_*(r))
 \end{equation*}
holds.

{\bf Theorem 1.} {\it Let $f,T, \Lambda$ satisfy  conditions {\rm (\ref{CondLambda})--(\ref{Condomega})}, and $a = \|f(x_0)\|$ satisfies the inequality}
 \begin{equation}\label{Cond for a}
 \lambda(r)\theta(r)w(r,a) \le r.
 \end{equation}

{\it Then equation {\rm (\ref{equation})} has a solution $x_* \in B(x_0,R)$, approximations {\rm (\ref{approximations})} converge to this solution $x_*$,  and the following inequalities  hold
 \begin{equation}\label{aposteriori}
 \|x_n - x_*\| \le \lambda(r)\theta(r)w(r,\|f(x_n)\|) \qquad (n =  0,1,2,\ldots),
 \end{equation}
({\sl aposteriori} estimates) and
 \begin{equation}\label{apriori}
 \|x_n - x_*\| \le \lambda(r)\theta(r)w(r,d^{(n)}(r,a)) \qquad (n =  0,1,2,\ldots),
 \end{equation}
({\sl apriori} estimates); moreover,}
 \begin{equation}\label{VeloLexp}
 {\rm Velo} \, (x_n)_{n=0}^\infty \le \frac{d(r,a)}{a}, \qquad  {\rm Lexp} \, (x_n)_{n=0}^\infty \le \mu'(0).
 \end{equation}

Here ${\rm Velo}$ and ${\rm Lexp}$ are the velocity of the rate of convergence and Lyapunov's exponent for the sequence $(x_n)_{n=0}^\infty$:
 \begin{equation*}
 {\rm Velo} \, (x_n)_{n=0}^\infty = \sup_n \, \frac{\|x_{n+1} - x_*\|}{\|x_n - x_*\|}, \qquad {\rm Lexp} \, (x_n)_{n=0}^\infty = \limsup_{n \to \infty} \, \sqrt[n]{\|x_n - x_*\|}.
 \end{equation*}

{\bf Proof.} Let us suppose first that  $r$, $0 < r \le \psi_*(r)$, is a number such that
  \begin{equation}\label{ineq with r}
  \|x_n - x_0\| \le r \qquad (n = 0,1,2,\ldots).
  \end{equation}

In this case, we remark, that equation (\ref{approximations}) implies the following identity
 \begin{equation}\label{BasicIdentity}
 \begin{array}{c} f(x_{n+1}) = f(x_n) - \Lambda(x_n,f(x_n))f'(x_n)T(x_n)f(x_n) + \\[10pt]  f(x_{n+1}) - f(x_n) - f'(x_n)(x_{n+1} - x_n) \qquad (n = 0,1,2,\ldots).\end{array}
 \end{equation}
This equation implies, in turn, the inequality
 \begin{equation}\label{BasicInequality}
 \begin{array}{c} \|f(x_{n+1})\| \le \|f(x_n) - \Lambda(x_n,f(x_n))f'(x_n)T(x_n)f(x_n)\| + \\[10pt]  \|f(x_{n+1}) - f(x_n) - f'(x_n)(x_{n+1} - x_n)\| \qquad (n = 0,1,2,\ldots).\end{array}
 \end{equation}

We can estimate the first term by means of condition (\ref{CondLambdaT}) ($h = f(x_n)$):
 \begin{equation}\label{LinearPart}
 \|f(x_n) - \Lambda(x_n,f(x_n))f'(x_n)T(x_n)f(x_n)\| \le \mu(r)\|f(x_n)\| \qquad (n = 0,1,2,\ldots).
 \end{equation}

In order to estimate the second term in (\ref{BasicInequality}) we note  first that
 \begin{equation*}
 \|f(\widetilde{x}) - f(x) - f'(x)(\widetilde{x} - x)\| \le \int\limits_0^1 \|f'((1 - t)x + t\widetilde{x}) - f'(x)\| \, \|\widetilde{x} - x\| \, dt
 \end{equation*}
and further, due to (\ref{Condomega}), one has
 \begin{equation*}
 \|f(\widetilde{x}) - f(x) - f'(x)(\widetilde{x} - x)\| \le \int\limits_0^1 \omega(r,t(\|\widetilde{x} - x\|)) \, \|\widetilde{x} - x\|\, dt =
 \end{equation*}
 \begin{equation*}
 \int\limits_0^{\|\widetilde{x} - x\|} \omega(r,\tau) \, d\tau = \Omega(r,\|\widetilde{x} - x\|).
 \end{equation*}
From this inequality we obtain that
 \begin{equation*}
 \|f(x_{n+1}) - f(x_n) - f'(x_n)(x_{n+1} - x_n)\| \le \Omega(r,\|x_{n+1} - x_n\|) \qquad (n = 0,1,2,\ldots)
 \end{equation*}
and, further (\ref{CondLambda}) implies
 \begin{equation}\label{NonlinearPart}
 \|f(x_{n+1}) - f(x_n) - f'(x_n)(x_{n+1} - x_n)\| \le \Omega(r,\lambda(r)\theta(r)\|f(x_n)\|) \qquad (n = 0,1,2,\ldots).
 \end{equation}

As result, from (\ref{LinearPart}) and (\ref{NonlinearPart}) and the definition of $d(r,\phi)$ (see (\ref{Fun d})) we have the following estimate
 \begin{equation}\label{Relax}
 \|f(x_{n+1})\| \le d(r,\|f(x_n)\|) \qquad (n = 0,1,2,\ldots)
 \end{equation}
and
 \begin{equation}\label{Correction}
 \|x_{n+1} - x_n\| \le \lambda(r)\theta(r)d(r,\|f(x_n)\|) \qquad (n = 0,1,2,\ldots).
 \end{equation}

Recall that all inequalities (\ref{Relax})--(\ref{Correction}) were obtained under {\small\sc apriori} assumption (\ref{ineq with r}).

Inequality (\ref{ineq with r}) for $n = 1$ is evident:
 \begin{equation*}
 \|x_1 - x_0\| = \|\Lambda(x_0,f(x_0))T(x_0)f(x_0)\| \le \lambda(r)\theta(r)a \le \lambda(r)\theta(r)w(r,a) \le r.
 \end{equation*}
In a similar way, due to (\ref{Relax}) for $n = 1$ one has
  \begin{equation*}
 \|x_2 - x_1\| = \|\Lambda(x_1,f(x_1))T(x_1)f(x_1)\| \le \lambda(r)\theta(r)\|f(x_1)\| \le
 \end{equation*}
 \begin{equation*}
 \lambda(r)\theta(r)d(r,\|f(x_0)\| \le \lambda(r)\theta(r)d(r,a),
 \end{equation*}
and, due to (\ref{Correction}) for $n = 2$ we obtain
 \begin{equation*}
 \|x_2 - x_0|| \le \|x_1 - x_0\| + \|x_2 - x_1\| \le  \lambda(r)\theta(r)(a + d(r,a)) \le \lambda(r)\theta(r)w(r,a) \le r.
 \end{equation*}
Repeating this  argument $n$ times we obtain
 \begin{equation}\label{Repeating}\begin{array}{c}
 \|x_{k+1} - x_k\| \le \|\Lambda(x_k,f(x_k))T(x_0)f(x_k)\| \le \lambda(r)\theta(r)\|f(x_k)\| \le \\[14pt] \lambda(r)\theta(r)d(r,\|f(x_{k-1})\|) \le \ldots \le \lambda(r)\theta(r)d^{(k)}(r,\|f(x_0)\|) = \lambda(r)\theta(r)d^{(k)}(r,a),
 \end{array}\end{equation}
and, due to (\ref{Correction}) for $k = 1,2,\ldots,n$ we have
 \begin{equation*}
 \|x_n - x_0\| \le \|x_1 - x_0\| + \|x_2 - x_1\| + \ldots + \|x_n - x_{n-1}\| \le
 \end{equation*}
 \begin{equation*}
 \lambda(r)\theta(r)(a + d(r,a) + \ldots + d^{(n-1)}(r,a)) \le \lambda(r)\theta(r)w(r,a) \le r.
 \end{equation*}

Thus, by induction, inequalities (\ref{ineq with r}) hold for all $n = 1,2,\ldots$.

Now let us verify  inequalities (\ref{aposteriori}) and (\ref{apriori}). For each $n = 0,1,2,\ldots$ due to (\ref{Relax}) (replacing   $n$ by  $k$) one obtains
 \begin{equation*}
 \|x_n - x_*\| \le \D\sum_{k=n}^\infty \|x_{k+1} - x_k\| \le \lambda(r)\theta(r) \D\sum_{k=n}^\infty d^{(k-n)}(r,\|f(x_n)\|) =
 \end{equation*}
 \begin{equation*}
 \lambda(r)\theta(r) \D\sum_{k=0}^\infty d^{(k)}(r,\|f(x_n)\|) = \lambda(r)\theta(r) \, w(r,\|f(x_n)\|).
 \end{equation*}
The latter means that (\ref{aposteriori}) is true. Iterating (\ref{Relax}) we obtain the inequality
 \begin{equation*}
 \|f(x_n)\| \le d^{(n)}(r,a) \qquad (n = 0,1,2,\ldots)
 \end{equation*}
and we  see that (\ref{apriori}) follows from (\ref{aposteriori}) and this inequality. Replacing (\ref{aposteriori}) by (\ref{Repeating}) we obtain  (\ref{apriori}).

Inequalities (\ref{VeloLexp}) follow from Figure presented above.

The proof  is complete.

\vspace{0.5cm}

{\bf 2.} In this section we discuss some problems concerning conditions (\ref{CondLambda}) and (\ref{CondLambdaT}). Since $X$ is a Hilbert space one can replace  condition (\ref{CondLambdaT}) by  the following inequality
 \begin{equation}\label{CondLambdaTmod}
 \begin{array}{c}\|f(x)\|^2 + 2\Lambda(x,f(x))(f(x),f'(x)T(x)f(x))  + \Lambda(x,f(x))^2\|f'(x)T(x)f(x)\|^2 \le \\[10pt] \mu^2(r)\|f(x)\|^2 \qquad (\|x - x_0\| \le r, \ 0 \le r \le R). \end{array}
 \end{equation}
Therefore, verification of (\ref{CondLambdaT}) is reduced to analysis of the following quadratic form
 \begin{equation}\label{QuadraticForm}
 \|f(x)\|^2 + 2\Lambda(f(x),f'(x)T(x)f(x))  + \Lambda^2\|f'(x)T(x)f(x)\|^2 \qquad (\|x - x_0\| \le R) \end{equation}
(with respect to $\Lambda = \Lambda(x,f(x))$). Moreover, we can see here in what way the choice of the functional $\Lambda(x,f(x))$  influences  the value of the coefficient $\mu(r)$. Naturally, we must try to choose  this functional in such a way  that the coefficient $\mu(r)$ would be  less than $1$ and, when possible, would be  minimal.

We remark that the inequality $\mu(r) < 1$ is possible only in the case when the discriminant of the polynomial
 \begin{equation*}
 (1 - \mu^2(r))\|f(x)\|^2 + 2\Lambda(f(x),f'(x)T(x)f(x))  + \Lambda^2\|f'(x)T(x)f(x)\|^2
 \end{equation*}
(with respect to $\Lambda$) is positive:
 \begin{equation*}
 (f(x),f'(x)T(x)f(x))^2 - (1 - \mu^2(r)) \ \|f(x)\|^2 \ \|f'(x)T(x)f(x)\|^2 > 0.
 \end{equation*}
This inequality implies that $(f(x),f'(x)T(x)f(x)) \ne 0$ for all $x$, $\|x - x_0\| \le R$ and for all $f(x)$, $f(x) \ne 0$. In what follows for definiteness sake we consider the case when $(f(x),f'(x)T(x)f(x)) > 0$ for all $x$, $\|x - x_0\| \le R$ with $f(x) \ne 0$. In other words, in what follows we assume that
 \begin{equation}\label{GenAss}
 (f(x),f'(x)T(x)f(x)) > 0 \qquad (\|x - x_0\| \le R, \ h \ne 0).
 \end{equation}

Let
 \begin{equation}\label{Defnu}
 \nu(r) = \inf_{\|x - x_0\| \le r} \ \frac{(f(x),f'(x)T(x)f(x))}{\|f(x)\| \, \|f'(x)T(x)f(x)\|} \qquad (0 < r \le R).
 \end{equation}
By virtue of (\ref{GenAss}) this function is nonnegative for all $0 \le r \le R$. It is evident that $\nu(r)$ is nonincreasing and $\nu(r) \le 1$. Hereafter  we assume that $0 < \nu(r) < 1$.

We consider two basic examples of the choice of  functional $\Lambda(x,h)$.

\emph{The first example}. \ The functional
 \begin{equation}\label{LambdaMin}
 \Lambda(x,h) = \frac{(h,f'(x)T(x)h)}{\|f'(x)T(x)h\|^2}
 \end{equation}
 makes the values of quadratic form (\ref{QuadraticForm}) minimal. These values are
 \begin{equation*}
 \|f(x)\|^2 - \frac{(f(x),f'(x)T(x)f(x))^2}{\|f'(x)T(x)f(x)\|^2}.
 \end{equation*}
Under this choice of $\Lambda(x,h)$ condition (\ref{CondLambdaT}) holds with
 \begin{equation}\label{MuMin}
 \mu(r) = \sqrt{1 - \nu^2(r)}.
 \end{equation}
In addition to this we have the equation
 \begin{equation}\label{DefTpKK}
 \lambda(r) = \sup_{\|x - x_0\| \le r} \ \frac{(f(x),f'(x)T(x)f(x))}{\|f'(x)T(x)f(x)\|^2}. \end{equation}

\emph{The second example}. \ The functional
 \begin{equation}\label{LambdaAl}
 \Lambda(x,h) = \frac{\|h\|^2}{\vartheta(h,f'(x)T(x)h)}
 \end{equation}
($\vartheta$ is a scalar parameter, $0 < \vartheta \le 2$) converts quadratic form in (\ref{QuadraticForm}) into the following expression
 \begin{equation*}
 \bigg(1 - \frac2\vartheta + \frac{\|f(x)\|^2\|f'(x)T(x)f(x)\|^2}{\vartheta^2(f(x),f'(x)T(x)f(x))^2}\bigg) \, \|f(x)\|^2.
 \end{equation*}
In this case, condition (\ref{CondLambdaT}) holds with
 \begin{equation}\label{MuAl}
 \mu(r) = \sqrt{1 - \frac2\vartheta + \frac1{\vartheta^2\nu^2(r)}}.
 \end{equation}
We note in this case that inequality $\mu(r) < 1$ is true only when
 \begin{equation}\label{nurest}
 \nu(r) > \sqrt{\frac1{2\vartheta}}, \qquad 0 \le r \le R.
 \end{equation}
Now, in this case, we have the equation
 \begin{equation}\label{DefTpKFA}
 \lambda(r) = \frac1\vartheta \, \sup_{\|x - x_0\| \le r} \ \frac{\|f(x)\|^2}{(f(x),f'(x)T(x)f(x))}.
 \end{equation}

Thus, the verification of condition (\ref{CondLambdaT}) in basic cases is reduced to the calculation of function (\ref{Defnu}) and defining its bounds. Remark that $\nu(r) \ge \widetilde{\nu}(r)$, where
 \begin{equation}\label{DefTnu}
 \widetilde{\nu}(r) = \inf_{\|x - x_0\| \le r, \ \|h\| = 1} \ \frac{(h,f'(x)T(x)h)}{\|h\| \, \|f'(x)T(x)h\|} \qquad (0 < r \le R).
 \end{equation}
Similarly, the verification of condition (\ref{CondLambdaT}) in basic cases is reduced to the analysis of functions (\ref{DefTpKK}) and (\ref{DefTpKFA}). Remark here that $\lambda(r) \le \widetilde{\lambda}(r)$, where
 \begin{equation}\label{DefTpiKK}
 \widetilde{\lambda}(r) = \sup_{\|x - x_0\| \le r, \ \|h\| = 1} \ \frac{(h,f'(x)T(x)h)}{\|f'(x)T(x)h\|^2}
 \end{equation}
in the case when $\Lambda(x,h)$ is defined by  (\ref{LambdaMin}) and
 \begin{equation}\label{DefTpiKFA}
 \widetilde{\lambda}(r) = \frac1\vartheta \, \sup_{\|x - x_0\| \le r, \ \|h\| = 1} \ \frac{\|h\|^2}{(h,f'(x)T(x)h)}
 \end{equation}
in the case when $\Lambda(x,h)$ is defined by (\ref{LambdaAl}).

Let us clarify the sense of magnitudes defined by  formulas (\ref{DefTnu}), (\ref{DefTpiKK}), and (\ref{DefTpiKFA}).

Inequalities $\widetilde{\nu}(r) > 0$, $0 < r \le R$, (see (\ref{DefTnu})) mean that operators $f'(x)T(x)$, \linebreak 
$\|x - x_0\| \le R$, belong to the class of operators $B: \ H \to H$ possessing  the property
 \begin{equation}\label{Acuteness}
 (h,Bh) \ge \nu \|h\| \, \|Bh\|, \qquad h \in H,
 \end{equation}
where $\nu$ is a positive constant, characterizing each operator $B$. As far as we know, this class of operators has not been  studied. Property (\ref{Acuteness}) for a fixed operator $B$ follows from the usual positive definiteness of the operator $B$. In its turn, it implies the property of positivity in the Krasnosel'ski\u{\i} -- Samarski\u{\i} sense. Geometrically inequality (\ref{Acuteness}) means that the angle between vectors $h$ and $Bh$ is acute.

Inequalities $\widetilde{\lambda}(r) < \infty$, $0 < r \le R$, with $\widetilde{\lambda}(r)$ defined by   (\ref{DefTpiKK}), mean that operators $f'(x)T(x)$, $\|x - x_0\| \le R$, belong to the class of operators $B: \ H \to H$ with the property
 \begin{equation}\label{Positivity??}
 (h,Bh) \le \lambda \|Bh\|^2, \qquad h \in H,
 \end{equation}
where $\lambda$ is a positive constant, characterizing each operator $B$. This class of operators
has not been  studied as well. Remark that operators satisfying both properties (\ref{Acuteness}) and (\ref{Positivity??}) have closed range; this situation holds, if functional $\Lambda(x,h)$ is defined by  formula (\ref{LambdaMin}).

Inequalities $\widetilde{\lambda}(r) < \infty$, $0 \le r \le R$, with $\widetilde{\lambda}(r)$ defined by  (\ref{DefTpiKFA}), mean that operators $f'(x)T(x)$, $\|x - x_0\| \le R$, are strongly positive definite in the standard sense (however, these operators need not necessarily  be self-adjoint).

\vspace{0.5cm}

{\bf 3.} Theorem 1 can be considered as a general principle of convergence for methods of gradient type. In particular, it includes numerous results for classical variants of method of minimal residuals, method of steepest descent, method of minimal errors and some others. Below, we present some of them.

{\bf Example 1: Method of minimal residuals}, see \cite{KK} and also \cite{KU,1,3}. In this case approximations are calculated by the formulas
 \begin{equation}\label{Mmr(1)}
 x_{n+1} = x_n - \frac{(f(x_n),f'(x_n)f(x_n))}{\|f'(x_n)f(x_n)\|^2} \ f(x_n) \quad (n =  1,2,\ldots).
 \end{equation}
Functional $\Lambda(x,h)$ correspondingly to them  and operators $T(x)$ are defined by means of the equations
 \begin{equation*}
 \Lambda(x,h) =  \frac{(h,f'(x)h)}{\|f'(x)h\|^2}, \qquad T(x) = I.
 \end{equation*}
Hence
 \begin{equation*}
 \lambda(r) = \sup_{\|x - x_0\| \le r} \ \frac{(f(x),f'(x)f(x))}{\|f'(x)f(x)\|^2}, \qquad \nu(r) = \inf_{\|x - x_0\| \le r} \ \frac{(f(x),f'(x)f(x))}{\|f(x)\| \, \|f'(x)f(x)\|}
 \end{equation*}
and
 \begin{equation}\label{munu1}
 \widetilde{\lambda}(r) = \sup_{\|x - x_0\| \le r}  \sup_{\|h\| = 1} \ \frac{(h,f'(x)h)}{\|f'(x)h\|^2}, \qquad \widetilde{\nu}(r) = \inf_{\|x - x_0\| \le r}  \inf_{\|h\| = 1} \ \frac{(h,f'(x)h)}{\|h\| \, \|f'(x)h\|}.
 \end{equation}

{\bf Example 2: Method of minimal co-errors}, see \cite{6}. In this case approximations are defined by the formulas
 \begin{equation}\label{Mmr(2)}
 x_{n+1} = x_n - \frac{\|{f'}^*(x_n)f(x_n)\|^2}{\|f'(x_n){f'}^*(x_n)f(x_n)\|^2} \ {f'}^*(x_n)f(x_n) \quad (n = 1,2,\ldots).
 \end{equation}
Functional $\Lambda(x,h)$ correspondingly to them  and operators $T(x)$ are defined by means of the equations
 \begin{equation*}
 \Lambda(x,h) =  \frac{\|{f'}^*(x)h\|^2}{\|f'(x){f'}^*(x)h\|^2}, \qquad T(x) = {f'}^*(x).
 \end{equation*}
Hence
 \begin{equation*}
 \lambda(r) = \sup_{\|x - x_0\| \le r} \D\frac{\|{f'}^*(x)f(x)\|^2}{\|f'(x){f'}^*(x)f(x)\|^2}, \quad  \nu(r) = \inf_{\|x - x_0\| \le r}  \ \frac{\|{f'}^*(x)f(x)\|^2}{\|f(x)\| \, \|f'(x){f'}^*(x)f(x)\|}
 \end{equation*}
and
 \begin{equation}\label{munu2}
 \widetilde{\lambda}(r) = \sup_{\|x - x_0\| \le r}  \sup_{\|h\| = 1} \ \frac{\|{f'}^*(x)h\|^2}{\|f'(x){f'}^*(x)h\|^2}, \quad \widetilde{\nu}(r) = \inf_{\|x - x_0\| \le r}  \inf_{\|h\| = 1} \ \frac{\|{f'}^*(x)h\|^2}{\|h\| \, \|f'(x){f'}^*(x)h\|}.
 \end{equation}

{\bf Example 3: Method of steepest descent}, see \cite{KA}, also \cite{FR2,KU,2,5}. In this case approximations are defined by the formulas
 \begin{equation}\label{Msd}
 x_{n+1} = x_n - \frac{\|f(x_n)\|^2}{(f(x_n),f'(x_n)f(x_n))} \ f(x_n) \quad (n =  1,2,\ldots)
 \end{equation}
Functional $\Lambda(x,h)$ corresponding to them  and operators $T(x)$ are defined by means of the equations
 \begin{equation*}
 \Lambda(x,h) = \frac{\|h\|^2}{(h,f'(x)h)}, \qquad T(x) = I.
 \end{equation*}
Hence
 \begin{equation*}
 \lambda(r) = \sup_{\|x - x_0\| \le r} \D\frac{\|f(x)\|^2}{(f(x),f'(x)f(x))}, \qquad \nu(r) = \inf_{\|x - x_0\| \le r}  \ \frac{(f(x),f'(x)f(x))}{\|f(x)\| \, \|f'(x)f(x)\|}
 \end{equation*}
and
 \begin{equation}\label{munu3}
 \widetilde{\lambda}(r) = \sup_{\|x - x_0\| \le r}  \sup_{\|h\| = 1} \ \frac{\|h\|^2}{(h,f'(x)h)}, \qquad \widetilde{\nu}(r) = \inf_{\|x - x_0\| \le r}  \inf_{\|h\| = 1} \ \frac{(h,f'(x)h)}{\|h\| \, \|f'(x)h\|}.
 \end{equation}

{\bf Example 4: Altman's Method of steepest descent,} see \cite{9,10,11}. In this case approximations are defined by the formulas
 \begin{equation}\label{MsdA}
 x_{n+1} = x_n - \frac{\|f(x_n)\|^2}{\vartheta(f(x_n),f'(x_n)f(x_n))} \ f(x_n)  \quad (n = 1,2,\ldots)
 \end{equation}
Functional $\Lambda(x,h)$ corresponding to them  and operators $T(x)$ are defined by means of the equations
 \begin{equation*}
 \Lambda(x,h) = \frac{\|h\|^2}{\vartheta(h,f'(x)h)}, \qquad T(x) = I.
 \end{equation*}
Hence
 \begin{equation*}
 \lambda(r) = \sup_{\|x - x_0\| \le r} \D\frac{\|f(x)\|^2}{\vartheta(f(x),f'(x)f(x))}, \qquad \nu(r) = \inf_{\|x - x_0\| \le r}  \ \frac{(f(x),f'(x)f(x))}{\|f(x)\| \, \|f'(x)f(x)\|}
 \end{equation*}
and
 \begin{equation}\label{munu4}
 \widetilde{\lambda}(r) = \sup_{\|x - x_0\| \le r}  \sup_{\|h\| = 1} \ \frac{\|h\|^2}{\vartheta(h,f'(x)h)}, \qquad \widetilde{\nu}(r) = \inf_{\|x - x_0\| \le r}  \inf_{\|h\| = 1} \ \frac{(h,f'(x)h)}{\|h\| \, \|f'(x)h\|}.
 \end{equation}

{\bf Example 5: Method of minimal errors}, see \cite{FR1,FR3,FR4}. In this case approximations are defined by
the formulas
 \begin{equation}\label{Mme}
 x_{n+1} = x_n - \frac{\|f(x_n)\|^2}{\|{f'}^*(x_n)f(x_n)\|^2} \ {f'}^*(x_n)f(x_n)  \quad (n = 1,2,\ldots)
 \end{equation}
Functional $\Lambda(x,h)$ corresponding to them  and operators $T(x)$ are defined by means of the equations
 \begin{equation*}
 \Lambda(x,h) = \frac{\|h\|^2}{\|{f'}^*(x)h\|^2}, \qquad T(x) = {f'}^*(x)
 \end{equation*}
Hence
 \begin{equation*}
 \lambda(r) = \sup_{\|x - x_0\| \le r} \frac{\|f(x)\|^2}{\|{f'}^*(x)f(x)\|^2}, \qquad \nu(r) = \sup_{\|x - x_0\| \le r}  \ \frac{\|{f'}^*(x)f(x)\|^2}{\|f(x)\| \, \|f'(x){f'}^*(x)f(x)\|}
 \end{equation*}
and
 \begin{equation}\label{munu5}
 \widetilde{\lambda}(r) = \sup_{\|x - x_0\| \le r}  \sup_{\|h\| = 1} \ \frac{\|h\|^2}{\|{f'}^*(x)h\|^2}, \qquad \widetilde{\nu}(r) = \inf_{\|x - x_0\| \le r}  \inf_{\|h\| = 1} \ \frac{\|{f'}^*(x)h\|^2}{\|h\| \, \|f'(x){f'}^*(x)h\|}.
 \end{equation}

{\bf Example 6: Altman's Method of minimal errors}, see \cite{9,10,11}. In this case approximations are defined by formulas
 \begin{equation}\label{MmeA}
 x_{n+1} = x_n - \frac{\|f(x_n)\|^2}{\vartheta\|{f'}^*(x_n)f(x_n)\|^2} \ {f'}^*(x_n)f(x_n)  \quad (n = 1,2,\ldots)
 \end{equation}
Functional $\Lambda(x,h)$ correspondingly to them  and operators $T(x)$ are defined by means of the equations
 \begin{equation*}
 \Lambda(x,h) = \frac{\|h\|^2}{\vartheta\|{f'}^*(x)h\|^2}, \qquad T(x) = {f'}^*(x)
 \end{equation*}
Hence
 \begin{equation*}
 \lambda(r) = \sup_{\|x - x_0\| \le r} \frac{\|f(x)\|^2}{\vartheta\|{f'}^*(x)f(x)\|^2}, \qquad \nu(r) = \sup_{\|x - x_0\| \le r}  \ \frac{\|{f'}^*(x)f(x)\|^2}{\|f(x)\| \, \|f'(x){f'}^*(x)f(x)\|}
 \end{equation*}
and
 \begin{equation}\label{munu6}
 \widetilde{\lambda}(r) = \sup_{\|x - x_0\| \le r}  \sup_{\|h\| = 1} \ \frac{\|h\|^2}{\vartheta\|{f'}^*(x)h\|^2}, \qquad \widetilde{\nu}(r) = \inf_{\|x - x_0\| \le r}  \inf_{\|h\| = 1} \ \frac{\|{f'}^*(x)h\|^2}{\|h\| \, \|f'(x){f'}^*(x)h\|}.
 \end{equation}

\vspace{0.5cm}

{\bf 4.} It is natural to study possibilities of generalizations of the argument in {\bf 1--3}
which is  specific for Hilbert spaces onto the case of operator equations on Banach spaces. As was pointed out in {\bf 3} the assumption that $H$ is a Hilbert space was used for the passage from condition (\ref{CondLambdaT}) to the analysis of the quadratic form in (\ref{CondLambdaTmod}). This passage is based on the classical identity for Hilbert spaces
 \begin{equation}\label{SquareEq}
 \|x + y\|^2 = \|x\|^2 + 2(x,y) + \|y\|^2.
 \end{equation}
It is known that this identity is not true for Banach spaces. Moreover, this identity contains the scalar product and this means that all constructions of {\bf 1--3} are impossible in Banach space.

However, we can try to use Lumer's semiscalar product $[x,y]$ that can be defined on an arbitrary Banach space (see \cite{LUM}). Let us recall the corresponding definitions.

Let $X$ be a Banach space, $X^*$ its conjugate. The set
 \begin{equation*}
 {\frak J}x = \{l \in X^*: \ \|l\| = \|x\|, \ (l,x) =  \|x\|^2\}
 \end{equation*}
is nonempty due to Hahn -- Banach theorem. If we have a selection
 \begin{equation*}
 J: \ X \to X^*, \quad Jx \in {\frak J}x
 \end{equation*}
of multifunction ${\frak J}: \ X \mapsto X^*$, we can define the {\it semiscalar product}
 \begin{equation}\label{LumProd}
 [x,y] = \langle Jx,y \rangle.
 \end{equation}
The following properties of the product $[x,y]: \ X \times X \to {\Bbb R}$ are true:
\begin{enumerate}\parskip = -2pt
\item[(a)] $[x,x] = \|x\|^2$,
\item[(b)] $[\lambda x,y] = \lambda [x,y]$,
\item[(c)] $[x,\alpha_1 y_1 + \alpha_2 y_2] = \alpha_1[x,y_1] + \alpha_2[x,y_2]$,
\item[(d)] $[x,y] \le \|x\|\|y\|.$
\end{enumerate}

In articles \cite{T,TP,YU,CIOR}  the equivalence of the following properties of a Banach space  $X$  is proved:
\begin{enumerate}\parskip = -2pt
\item[(i)] $\|x + y\|^2 + \sigma\|x - y\|^2 \ge (\|x\|^2 + \|y\|^2) \ \ (x, y \in X)$,
\item[(ii)] $\|x + y\|^2 + \|x - y\|^2 \le 2(\|x\| + \sigma\|y\|^2) \ \ (x, y \in X)$,
\item[(iii)] $\|(1 - \lambda)x + \lambda y\|^2 \ge (1 - \lambda)\|x\|^2 + \lambda\|y\|^2 - \sigma(1 - \lambda)\lambda\|x - y\|^2  \ \ (x, y \in X, \ 0 < \lambda < 1),$
\item[(iv)] $\|x + y\|^2 \le \|x\|^2 + 2[x,y] + \sigma\|y\|^2  \ \ (x, y \in X)$,
\item[(v)] $\langle Jx - Jy,x - y \rangle \le \sigma\|x - y\|^2 \ \ (x, y \in X)$,
\item[(vi)] $\|Jx - Jy\|  \le \sigma\|x - y\| \ \ (x, y \in X);$
\end{enumerate}
here $\sigma$ is a number, $\sigma \ge 1$. Banach spaces with these properties are known as the  {\it spaces with Bynum's property}. The spaces $L_p$ ($2 \le p < \infty$), some Orlich spaces, and spaces with Hanner inequality are spaces with Bynum's property. Remark that in the case of the spaces with Bynum property the multifunction ${\frak J}$ is a usual single-valued function, or, in other words, in these  spaces the semiscalar product is defined uniquely.

Thus, if $X$ is a Banach space with Bynum property, we can use the inequality
 \begin{equation}\label{SquareIneq}
 \|x + y\|^2 \le \|x\|^2 + 2[x,y] + \sigma\|y\|^2.
 \end{equation}
It is evident that this inequality can be used instead of  equality (\ref{SquareEq}) in all the reasoning and constructions presented above for the case of Hilbert space.

Now we return to operator equation (\ref{equation}), however, in the case when the left-hand side of this equation is an operator $f(x)$ on a Banach space $X$ with the Bynum property. The semiscalar product $[x,y]$ in this space allows us to consider instead of inequality (\ref{CondLambdaTmod}) the analogous inequality
 \begin{equation}\label{BQuadFormIneq}
 \begin{array}{c}\|h\|^2 + 2\Lambda(x,h)[h,f'(x)T(x)h]  + \sigma\Lambda(x,h)^2\|f'(x)T(x)h\|^2 \le \mu^2(r)\|h\|^2 \\[10pt] (0 \le r \le r, \ 0 \le r \le R).\end{array}
 \end{equation}
Now we can repeat the argument of the previous sections. In particular, we must investigate the left-hand side of this inequality:
 \begin{equation}\label{QuadraticFormBS}
 \|f(x)\|^2 + 2\Lambda[f(x),f'(x)T(x)f(x)] + \sigma\Lambda^2\|f'(x)T(x)f(x)\|^2 \qquad (\|x - x_0\| \le R).
 \end{equation}

Below we formulate an analogue to Theorem 1 for operator equations on Banach space $X$. Naturally, we deal with approximations (\ref{approximations}) under the same assumptions (\ref{CondLambda})--(\ref{Condomega}). However, instead of function (\ref{Fun d}) we must consider the function
 \begin{equation}\label{6!}
 d_\sigma(r,\phi) = \mu(r)\phi + \sigma\Omega(r,\lambda(r)\theta(r)\phi) \qquad \biggl(\Omega(r,t) = \int\limits_0^t \omega(r,\tau) \, d\tau\biggr);
 \end{equation}
functions $d^{(n)}(r,\phi)$, $n = 1,2,\ldots$, and $w(r,\phi)$ are also replaced by   the corresponding functions $d_\sigma^{(n)}(r,\phi)$, $n = 1,2,\ldots$, and $w_\sigma(r,\phi)$. In the proof of this analogue, only formula (18) is changed essentially; of course functions $d^{(n)}(r,\phi)$, $n = 1,2,\ldots$ are also replaced by  functions $d_\sigma^{(n)}(r,\phi)$, $n = 1,2,\ldots$, and $w_\sigma(r,\phi)$.

{\bf Theorem 2.} {\it Let $f,T, \Lambda$ satisfy  conditions {\rm (\ref{CondLambda})--(\ref{Condomega})}, and $a = \|f(x_0)\|$ satisfies the inequality}
 \begin{equation}\label{7!}
 \lambda(r)\theta(r)w_\sigma(r,a) \le r.
 \end{equation}

{\it Then equation {\rm (\ref{equation})} has a solution $x_* \in B(x_0,R)$, approximations {\rm (\ref{approximations})} converge to this solution $x_*$,  and the following inequalities  hold
 \begin{equation}\label{8!}
 \|x_n - x_*\| \le \lambda(r)\theta(r)w_\sigma(r,\|f(x_n)\|) \qquad (n =  0,1,2,\ldots),
 \end{equation}
({\sl aposteriori} estimates) and
 \begin{equation}\label{9!}
 \|x_n - x_*\| \le \lambda(r)\theta(r)w_\sigma(r,d_\sigma^{(n)}(r,a)) \qquad (n =  0,1,2,\ldots),
 \end{equation}
({\sl apriori} estimates); moreover,}
 \begin{equation}\label{10!}
 {\rm Velo} \, (x_n)_{n=0}^\infty \le \frac{d_\sigma(r,a)}{a}, \qquad  {\rm Lexp} \, (x_n)_{n=0}^\infty \le \mu'(0).
 \end{equation}

Now we can write natural analogues for concrete iterative methods from Section {\bf 3}. First, we remark that the analogues to(\ref{LambdaMin}) and (\ref{LambdaAl}) are
 \begin{equation}\label{LambdaMinB}
 \Lambda(x,h) = \frac{[h,f'(x)T(x)h]}{\sigma\|f'(x)T(x)h\|^2}
 \end{equation}
and
 \begin{equation}\label{LambdaAlB}
 \Lambda(x,h) = \frac{\|h\|^2}{\vartheta[h,f'(x)T(x)h]}.
 \end{equation}
So, we obtain the following analogue to  method (\ref{Mmr(1)}) in {\bf Example 1} (method of minimal residuals)
 \begin{equation}\label{3!}
 x_{n+1} = x_n - \frac{[f(x_n),f'(x_n)f(x_n)]}{\sigma\|f'(x_n)f(x_n)\|^2} f(x_n)  \qquad (n = 0,1,2,\ldots),
 \end{equation}
and analogues to method (\ref{Msd}) in {\bf Example 3} (method of steepest descent)
 \begin{equation}\label{4!}
 x_{n+1} = x_n - \frac{\|f'(x_n)f(x_n)\|^2}{[f(x_n),f'(x_n)f(x_n)]} f(x_n) \qquad (n = 0,1,2,\ldots),
 \end{equation}
and to method (\ref{MsdA}) in {\bf Example 4} (Altman method of steepest descent)
 \begin{equation}\label{5!}
 x_{n+1} = x_n - \frac{\|f'(x_n)f(x_n)\|^2}{\vartheta[f(x_n),f'(x_n)f(x_n)]} f(x_n) \qquad (n = 0,1,2,\ldots).
 \end{equation}
In all these cases $T(x) = I$ for all $x$, $\|x - x_0\| \le R$.

Equations (\ref{munu1}) for method (\ref{Mmr(1)}) must be replaced by the  equations
 \begin{equation}\label{munu1B}
 \widetilde{\lambda}(r) = \sup_{\|x - x_0\| \le r} \sup_{\|h\| = 1} \ \frac{[h,f'(x)h]}{\vartheta\|f'(x)h)\|^2}, \qquad \widetilde{\nu}(r) = \inf_{\|x - x_0\| \le r} \inf_{\|h\| = 1} \ \frac{[h,f'(x)h]}{\|f'(x)h\|^2}.
 \end{equation}
Analogously, equations (\ref{munu3}) for method (\ref{Msd}) must be replaced by the  equations
 \begin{equation}\label{13!}
 \widetilde{\lambda}(r) = \sup_{\|x - x_0\| \le r} \sup_{\|h\| = 1} \ \frac{\|h\|^2}{[h,f'(x)h]}, \qquad \widetilde{\nu}(r) = \inf_{\|x - x_0\| \le r}  \inf_{\|h\| = 1} \ \frac{[h,f'(x)h]}{\|h\| \, \|f'(x)h\|},
 \end{equation}
and equations (\ref{munu4}) for method (\ref{MsdA}) must be replaced by
 \begin{equation}\label{15!}
 \widetilde{\lambda}(r) = \sup_{\|x - x_0\| \le r} \sup_{\|h\| = 1} \  \frac{\|h\|^2}{\vartheta[h,f'(x)h]}, \qquad \widetilde{\nu}(r) = \inf_{\|x - x_0\| \le r}  \inf_{\|h\| = 1} \ \frac{(h,f'(x)h)}{\|h\| \, \|f'(x)h\|}.
 \end{equation}

The approach put forward here is not applicable for methods in {\bf Examples 2}, {\bf 5}, and {\bf 6}, because the formal terms ${f'}^*(x)f(x)$, $\|x - x_0\| \le R$ are not operators on $X$, if $X$ is a Banach space.

The main idea of the passage from Hilbert to Banach space can be generalized on a more general class of Banach spaces. In particular, we can consider the class of Banach spaces where  the following inequality (see \cite{Xu})
 \begin{equation*}
 \|x + y\|^p \le \|x\|^p + p[x,y]_p + \sigma\|y\|^p
 \end{equation*}
holds, here $p \in (1,\infty)$, $\sigma \ge 1$,
 \begin{equation*}
 [x,y]_p = \langle J_px,y \rangle,
 \end{equation*}
and
 \begin{equation*}
 {\frak J}_px = \{l \in X^*: \ \|l\| = \|x\|^{p-1}, \ (l,x) =  \|x\|^p\},
 \end{equation*}
 \begin{equation*}
 J_p: \ X \to X^*, \quad J_px \in {\frak J}_px.
 \end{equation*}
The simple example of such spaces gives the spaces $L_p$, $1 < p \le 2$. We omit the formulations of the corresponding results.

\vspace{0.5cm}

\begin{center} {\bf Literature} \end{center} \vspace{-2.0cm}

\def\refname{}


\begin{thebibliography}{20}\parskip = -0.5pt

\bibitem{9}
{\bf Altman M.}: {\it On the approximate solution of non-linear functional equations}. Bull. Acad. polon. sci., Cl. 3, {\bf 5} (1957), \No 5, 457-460.

\bibitem{10}
{\bf Altman M.}: {\it An approximate method for solving linear equations in a Hilbert space}. Bull. Acad. polon. sci., Cl. 3, {\bf 5} (1957), \No 6, 601-604.

\bibitem{11}
{\bf Altman M.}: {\it On the approximate solutions of operator equations in Hilbert space}. Bull. Acad. polon. sci., Cl. 3, {\bf 5} (1957), \No 6, 605-609.

\bibitem{CIOR}
{\bf Cioranescu I.}: {\it Geometry of Banach spaces, duality mappings and nonlinear problems.} --- Kluwer Academic Publishers, Dordrecht - Boston - London, {\bf 62} (1990).

\bibitem{5}
{\bf Evkhuta O.N.}: {\it About methods of steepest descent for equations in Banach spaces}. Izvestia Vishish uchebnish zavedeni\u{\i} SKR, \No 11, (2005), 18-27

\bibitem{6}
{\bf Evkhuta O.N., Zabre\u{\i}ko P.P.}: {\it Some remarks on iterative gradient and quasi-gradient methods for approximate solution of nonlinear operator equations}. Izvestia Vishish uchebnish zavedeni\u{\i} SKR, 2005, \No 11, 27-37

\bibitem{FR1}
{\bf Fridman V.M.}: {\it New methods for solution of a linear operator equation}. --- Doklady AN SSSR, {\bf 128} (1959), \No 3, 482-484.

\bibitem{FR2}
{\bf Fridman V.M.}: {\it On the convergence of methods of steepest descent}. --- Uspekhi matem. nauk, {\bf 17} (1962), vyp. 3, 201-204.

\bibitem{FR3}
{\bf Fridman V.M.}: {\it An iterative process with minimal errors for a nonlinear operator equation}. --- Doklady AN SSSR, {\bf 139} \No 5 (1961), 1063-1066.

\bibitem{FR4}
{\bf Fridman V.M.}: {\it Method of minimal iterations with minimal errors for a system of linear algebraic equations with a symmetric matrix}. --- Zhurnal Vy\v{c}islit. Matem. i Matem. Fiziki, {\bf 2} (1962), \No 2, ?-?.

\bibitem{KA}
{\bf Kantorovi\v{c}, L.V. and Akilov. G.P.}: {\it  Functional Analysis}, Pergamon Press, New York, 1982.

\bibitem{2}
{\bf Kirsanova-Evkhuta O.N.}: {\it Convergence theorems of method of quickest descent and method of minimal errors}. Doklady NAN Belarusi, {\bf 48}, \No 2, (2004), 10-15

\bibitem{K1}
{\bf Kivistik L.V.}: {\it Оn the method of steepest descent for the solution of nonlinear equations.} --- Izvest. AN EstSSR, ser. fiz.-matem. i tekhn. nauk, {\bf 9} (1960), \No 2, 145-159.

\bibitem{K2}
{\bf Kivistik L.V.}: {\it On certain iterative methods for the solution of operator equations in Hilbert spaces.} --- Izvest. AN EstSSR, ser. fiz.-matem. i tekhn. nauk, {\bf 9} (1960), \No 3, 229-241.

\bibitem{K3}
{\bf Kivistik L.V.}: {\it On a modification of the iteration method with minimal residuals for the solution of nonlinear operator equations.} --- Doklady AN SSSR, {\bf 136} (1961), \No 1, 22-25.

\bibitem{KU}
{\bf Kivistik L.V., Ustaal A.}: {\it Some convergence theorems for iterative processes with minimal residuals}. --- U\v{c}en. Zapiski Tartussk. universiteta, 1962, \No 129, 382-393.

\bibitem{KK}
{\bf Krasnosel'ski\u{\i} M.A., Kre\u{\i}n S.G.}: {\it An iterative process with minimal residuals.} --- Matemati\v{c}eski\u{\i} sbornik, {\bf 31} (1952), \No 2, 315-334. (Russian)

\bibitem{KVZRS}
{\bf Krasnosel'ski\u{\i} M.A., Vainikko G.M., Zabreiko P.P., Rutitskii Y.B., Stetsenko V.Y.}: {\it Approximate solution of operator equations.} - Groningen, Wolters-Noordhoff Publishing, 1972. - 484 p.

\bibitem{LUM}
{\bf Lumer G.}: {\it Semi-inner-product spaces}. --- Transactions of the American Mathema-tical Society, {\bf 100} (1961), No. 1, 29-43

\bibitem{T}
{\bf Trubnikov Yu.V.} {\it Extremal Constructions in Nonsmooth Analysis and Operator Equations with Accretive Operators.} - Moscow: Astro-Press, 2002. - 256 p.

\bibitem{TP}
{\bf Trubnikov Yu.V., Perov A.I.}: {\it Differential Equations with Monotone Nonlinearities.} --- Minsk: Nauka i Tekhnika, 1986. - 200 с.

\bibitem{Xu}
{\bf Xu Hong-Kun}: {\it Inequalities in Banach Spaces with Applications.} --- Nonlinear Analysis: Theory, Methods, and Applications, {\bf 16} (1991), no 12, 1127-1138.

\bibitem{YU}
{\bf Yurgelas V.V.}: {\it On iterative processes of descent type and some properties of the dual mapping.} --- Voronezh, (1980), 1-30. Dep. VINITI 19.05.80, \No 1912--80c.

\bibitem{1}
{\bf Zabre\u{\i}ko P.P., Kirsanova-Evkhuta O.N.}: {\it A new theorem on convergence for the minimal residual method}. Vesti NAN Belarusi, ser. fiz.-matem. nauk, \No 2, (2004), 5-8.

\bibitem{3}
{\bf Zabre\u{\i}ko P.P., Kirsanova-Evkhuta O.N.}: {\it Method of minimal residuals in Banach spaces}. Doklady NAN Belarusi, {\bf 49}, \No 2, (2005), 5-10

\end{thebibliography}
\end{document}